# APPROXIMATING CONDITIONAL DISTRIBUTION FUNCTIONS USING DIMENSION REDUCTION

By Peter Hall and Qiwei Yao

*Australian National University and London School of Economics*

Motivated by applications to prediction and forecasting, we suggest methods for approximating the conditional distribution function of a random variable $Y$ given a dependent random $d$-vector $X$. The idea is to estimate not the distribution of $Y|X$, but that of $Y|\theta^{\mathrm{T}}X$, where the unit vector $\theta$ is selected so that the approximation is optimal under a least-squares criterion. We show that $\theta$ may be estimated root-$n$ consistently. Furthermore, estimation of the conditional distribution function of $Y$, given $\theta^{\mathrm{T}}X$, has the same first-order asymptotic properties that it would enjoy if $\theta$ were known. The proposed method is illustrated using both simulated and real-data examples, showing its effectiveness for both independent datasets and data from time series. Numerical work corroborates the theoretical result that $\theta$ can be estimated particularly accurately.

**1. Introduction.** Estimating a conditional distribution function is an important feature of many statistical problems, including, for example, regression analysis [see Yin and Cook (2002) and references therein], where a significant problem is prediction of a response for a given value of a multivariate explanatory variable. Specific applications include those in economics and finance [e.g., Foresi and Paracchi (1992), Bond and Patel (2000) and Watanabe (2000)], in signal processing and data mining [e.g., Adali, Liu and Sonmez (1997)] and a wide range of problems where forecasts are to be made from linear or nonlinear time-series [see, e.g., Chapter 10 of Fan and Yao (2003), and examples in Section 4 below].

In most of these applications one is interested in estimating the conditional distribution of a scalar random variable $Y$, given a random $d$-vector $X$. Even for small values of $d \geq 2$, a conventional nonparametric estimator can suffer poor accuracy, reflected in slow convergence rates. We suggest a solution to this difficulty, based on approximating the conditional distribution function of $Y$ given $X$ by that of $Y$ given $\theta^{\mathrm{T}}X$, where the unit $d$-vector $\theta$ is selected so that the approximation is optimal under an appropriate least-squares criterion. In particular, we avoid the problem of directly estimating the conditional distribution function of $Y$ given $X$.

Although we are dealing with a dimension-reduction problem, the object (i.e., the conditional distribution function) to be estimated is a function of both $\theta^{\mathrm{T}}x$ and $y$, while the index $\theta$ is a global parameter. This rules out the possibility of direct application of conventional dimension-reduction ideas, such as projection pursuit [e.g., Friedman and Stuetzle (1981), Friedman, Stuetzle and Schroeder (1984) and Huber (1985)] and single-index modeling techniques [e.g., Powell, Stock and Stoker (1989), Härdle, Hall and Ichimura (1993), Ichimura (1993) and Klein and Spady (1993)], which would lead to an estimator of $\theta$ depending on $y$. Instead we define a new criterion in terms of an accumulation of squared differences between the joint probabilities of $(Y, X)$ and the expected conditional probabilities of $Y$ given $\theta^{\mathrm{T}}X$, over a large class of subsets; see (2.2) and (2.4) in Section 2 below. Our search for the global parameter $\theta$ is based on leave-one-out local linear regression estimators for conditional distribution functions. Under very mild assumptions the resulting estimator $\hat{\theta}$ is root-$n$ consistent and asymptotically normally distributed.

Of course, our main purpose in computing $\hat{\theta}$ is so it can be used in a conditional distribution estimator. The root-$n$ convergence rate achieved by our estimator is so fast that the estimator of the conditional distribution function of $Y$, given $\hat{\theta}^{\mathrm{T}}X$, is first-order equivalent to its counterpart that would be used if the true value of $\theta$ were known.

The innovation and novelty of our methodology lie in the fact that we use dimension-reduction ideas to solve an important class of nonstandard multivariate nonparametric problems. We achieve this end by proposing new types of objective functions, with which are associated new theoretical and numerical properties. There exists an extensive literature on nonparametric estimation of conditional distributions. It includes work of Bhattacharya and Gangopadhyay (1990), Sheather and Marron (1990), Yu and Jones (1998) and Cai (2002) on conditional quantile regression; Rosenblatt (1969), Hyndman, Bashtannyk and Grunwald (1996), Fan, Yao and Tong (1996), Bashtannyk and Hyndman (2001) and Hyndman and Yao (2002) on conditional density estimation; and Hall, Wolff and Yao (1999) on estimation of conditional distribution functions. Dimension reduction has been discussed extensively in the context of regression and density approximation; in addition to the references cited earlier we mention the work of Friedman (1987), Jones and Sibson (1987), Li (1991) and Posse (1995).

This article is organized as follows. In Section 2 we introduce our method for estimating $\theta$. Asymptotic properties of estimators $\hat{\theta}$ and $\widehat{F}(\cdot|\theta^{\mathrm{T}}X)$ are presented in Section 3. Numerical examples involving both simulated models and a real-data application are given in Section 4. Technical arguments are outlined in Section 5.

## 2. Methodology.

2.1. *Motivation.* Assume we observe data $(X_i, Y_i)$, for $1 \le i \le n$, from the distribution of $(X, Y)$. Here $X$ is a $d$-vector and $Y$ is a scalar. Let $\Theta$ denote the set of $d$-variate unit vectors $\theta$ with first nonzero component positive, write $f$ for the density of $X$ and let $F_{Y|\theta^{\mathrm{T}}X}(\cdot|z)$ represent the distribution of $Y$ conditional on $\theta^{\mathrm{T}}X = z$. Given subsets $\mathcal{A}$ and $\mathcal{B}$ of $d$-dimensional space and of the real line, respectively, define

$$\pi_\theta(\mathcal{A}, \mathcal{B}) = \int_{\mathcal{A}} F_{Y|\theta^{\mathrm{T}}X}(\mathcal{B}|\theta^{\mathrm{T}}x)f(x)\,dx, \qquad \pi(\mathcal{A}, \mathcal{B}) = P(X \in \mathcal{A}, Y \in \mathcal{B}).$$

If, for some $\theta$ and all $x$, $F_{Y|\theta^{\mathrm{T}}X}(\cdot|\theta^{\mathrm{T}}x)$ is identical to the distribution of $Y$ given that $X = x$, then for this $\theta$, $\pi_\theta(\mathcal{A}, \mathcal{B}) = \pi(\mathcal{A}, \mathcal{B})$ for all $\mathcal{A}, \mathcal{B}$. We suggest taking the sets $\mathcal{A}$ to be $d$-variate spheres with differing centers and radii, and the sets $\mathcal{B}$ to be semi-infinite intervals.

We can estimate $F_{Y|\theta^{\mathrm{T}}X}$ using nonparametric methods, permitting us to estimate $\pi_\theta(\mathcal{A}, \mathcal{B})$. Of course, we can estimate $\pi(\mathcal{A}, \mathcal{B})$ as the proportion of pairs $(X_i, Y_i)$ that lie in $\mathcal{A} \times \mathcal{B}$. Hence, for each triple $(\theta, \mathcal{A}, \mathcal{B})$ we can estimate $\pi_\theta(\mathcal{A}, \mathcal{B})$ and $\pi(\mathcal{A}, \mathcal{B})$ under minimal conditions. (We shall denote estimators of $\pi_\theta$ and $\pi$ by $\hat{\pi}_\theta$ and $\hat{\pi}$, resp.) Therefore we can check (or, more formally, test) the hypothesis that $F_{Y|\theta^{\mathrm{T}}X}(\cdot|\theta^{\mathrm{T}}x)$ is identical to the distribution of $Y$ conditional on $X = x$, for all $x$, by examining the average value of $\{\hat{\pi}_\theta(\mathcal{A}, \mathcal{B}) - \hat{\pi}(\mathcal{A}, \mathcal{B})\}^2$ over a range of sets $\mathcal{A}$ and $\mathcal{B}$.

Although exact equality of $\pi$ and $\pi_\theta$ is unlikely in practice, the difference-based criterion noted above can be used to empirically select $\theta$ such that, in a global sense, the distribution of $Y$ given $\theta^{\mathrm{T}}X = \theta^{\mathrm{T}}x$ is a good approximation to the distribution of $Y$ given that $X = x$. Indeed, the argument in the previous paragraph suggests that methodology of this type could be based on the difference measure,

$$(2.1) \qquad S_1(\theta) = \iint \{\hat{\pi}_\theta(\mathcal{A}_\alpha, \mathcal{B}_\beta) - \hat{\pi}(\mathcal{A}_\alpha, \mathcal{B}_\beta)\}^2 w(\alpha, \beta)\,d\alpha\,d\beta,$$

where $w$ is a weight function and the integral is taken over a parameterization $(\alpha, \beta)$ of $(\mathcal{A}, \mathcal{B})$.

The spheres $\mathcal{A} = \mathcal{A}_\alpha$ should be such that the density $f_{\theta^{\mathrm{T}}X}$ of $\theta^{\mathrm{T}}X$ is bounded away from zero at all points $\theta^{\mathrm{T}}x$ with $\theta \in \Theta$ and $x \in \mathcal{A}_\alpha$. Otherwise, design sparseness problems can arise when nonparametrically estimating $F_{Y|\theta^{\mathrm{T}}X}$. Considerations of this type suggest taking the $\mathcal{A}_\alpha$'s to be $d$-variate spheres whose centers confine them to lie inside a larger, bounded region where $f$ is bounded away from zero. Such restrictions are unnecessary when considering the intervals $\mathcal{B}$, except that there is little point in giving emphasis to sets for which $P(Y \in \mathcal{B})$ is low.

For these reasons, when permitting $\mathcal{B}_\beta$ to be the interval $(-\infty, \beta)$ it is appropriate to take $w(\alpha, \beta)$ in (2.1) to be proportional to the density of $Y$ at $\beta$, and to not depend on $\alpha$. We shall achieve this end empirically, by replacing the double integral in (2.1) by a sum of integrals,

$$(2.2) \qquad S(\theta) = \sum_{j=1}^{n} \int \{\hat{\pi}_\theta(\mathcal{A}_\alpha, \mathcal{B}_{Y_j}) - \hat{\pi}(\mathcal{A}_\alpha, \mathcal{B}_{Y_j})\}^2 \, d\alpha,$$

where $\mathcal{B}_\beta$ denotes $(-\infty, \beta]$ and the integral is taken over an appropriate set of sphere centers and radii. In the future we shall use the notation $y$ instead of $\beta$.

When constructing our estimator of $F_{Y|\theta^{\mathrm{T}}X}$, which is central to our computation of $\hat{\pi}_\theta$, we shall use a "leave-one-out" technique, or more accurately, "leave-two-out." Our method will employ the empirical distribution of $\theta^{\mathrm{T}}X_i$ as a surrogate for the true distribution of $\theta^{\mathrm{T}}X$, and so, when $\theta^{\mathrm{T}}X_i$ appears in the argument of $\hat{F}_{Y|\theta^{\mathrm{T}}X}$, we shall omit $X_i$ from the latter. The second omission occurs because, as formula (2.2) suggests, we shall validate on $Y_j$ when constructing our least-squares criterion. Therefore we shall omit both the $i$th and the $j$th pairs when calculating $\hat{F}_{Y|\theta^{\mathrm{T}}X}$; see (2.3) below.

2.2. *Estimator of* $\theta$. With these principles in mind, let $h$ be a bandwidth and let $K$ be a kernel function, and define

$$T_{-i,-j}^{[k]}(\theta) = \frac{1}{(n-2)h} \sum_{i_1 : i_1 \neq i, j} K\left\{\frac{\theta^{\mathrm{T}}(X_i - X_{i_1})}{h}\right\} \left\{\frac{\theta^{\mathrm{T}}(X_i - X_{i_1})}{h}\right\}^k,$$

$$w_{i_1; -i, -j}(\theta) = K\left\{\frac{\theta^{\mathrm{T}}(X_i - X_{i_1})}{h}\right\}$$

$$(2.3) \qquad \qquad \times \left\{T_{-i,-j}^{[2]}(\theta) - \frac{\theta^{\mathrm{T}}(X_i - X_{i_1})}{h} T_{-i,-j}^{[1]}(\theta)\right\},$$

$$\hat{F}_{-i,-j}(y|\theta^{\mathrm{T}}X_i) = \left\{\sum_{i_1 : i_1 \neq i, j} w_{i_1; -i, -j}(\theta) I(Y_{i_1} \leq y)\right\}$$

$$\times \left\{\sum_{i_1 : i_1 \neq i, j} w_{i_1; -i, -j}(\theta)\right\}^{-1}.$$

Write simply $F(y|z)$ for $P(Y \leq y|\theta^{\mathrm{T}}X = z)$, and let $\mathcal{A}$ be a subset of $d$-variate space. In this notation, $\hat{F}_{-i,-j}(y|\theta^{\mathrm{T}}X_i)$ is a local linear estimator of $F(y|\theta^T X_i)$, based on data pairs other than the $i$th and the $j$th; and

$$\frac{1}{n-1} \sum_{i\,:\,i \neq j, X_i \in \mathcal{A}} \widehat{F}_{-i,-j}(y|\theta^T X_i)$$

is an estimator of $\pi_\theta(\mathcal{A}, \mathcal{B})$ when $\mathcal{B} = (-\infty, y]$.

As a rule we take $\alpha$ to be a $(d+1)$-vector, its first $d$ components denoting the center of $\mathcal{A}_\alpha$ and the last component, $r$ say, its radius. We suppose that $r \in \mathcal{J} = [r_1, r_2]$, where $0 \leq r_1 \leq r_2 \leq \infty$ and not both $r_1$ and $r_2$ vanish. In the case $r_1 = r_2$ the spheres all have the same radius, and here $\alpha$ should be interpreted as a $d$-vector, with integrals over $r$, in our discussion below, ignored. With this interpretation, our account of methodology applies to the case where $r$ takes values in the continuum as well as to that where $r$ is fixed. Clearly the latter instance can be generalized to the case of a finite number of discrete radii.

One approach is to average over all spheres $\mathcal{A}_\alpha$ that lie entirely within a given, fixed set $\mathcal{R}$. With this in mind, let $\mathcal{Q} = \{\alpha : \mathcal{A}_\alpha \subseteq \mathcal{R}\}$ be the set of sphere centers (and radii, if $\mathcal{J}$ is not degenerate). Write $\widehat{F}_{-j}(\mathcal{A}, y)$ for the proportion of the $n-1$ values of $(X_i, Y_i)$, for $i \neq j$, that satisfy $(X_i, Y_i) \in \mathcal{A} \times (-\infty, y]$. Put

$$
\begin{aligned}
S(\theta, \mathcal{A}) &= \sum_{j=1}^n \left\{ \widehat{F}_{-j}(\mathcal{A}, Y_j) - \frac{1}{n-1} \sum_{i\,:\,i \neq j, X_i \in \mathcal{A}} \widehat{F}_{-i,-j}(Y_j|\theta^T X_i) \right\}^2, \\
S(\theta) &= \int_{\mathcal{Q}} S(\theta, \mathcal{A}_\alpha)\, d\alpha.
\end{aligned}
\tag{2.4}
$$

The latter represents a particular form of $S(\theta)$ in (2.2). In practice, the integration over $\alpha$ in (2.4) is typically replaced by a sum over a class of selected balls; see (4.1) below. In fact the asymptotic theory in Section 3 still holds with this discrete version of $S(\theta)$ if the same replacement is applied wherever appropriate, including in condition (3.3).

We choose $\hat{\theta}$ to minimize $S(\theta)$ over $\theta \in \Theta$. Thus, $\hat{\theta}$ may be viewed as an estimator of $\theta_0$, the minimizer (over $\theta \in \Theta$) of

$$
S_0(\theta) = \int_{\mathcal{Q}} d\alpha \int \{F(\mathcal{A}_\alpha, y) - G_\theta(\mathcal{A}_\alpha, y)\}^2 f_Y(y)\, dy,
\tag{2.5}
$$

where $F(\mathcal{A}, y) = P\{(X, Y) \in \mathcal{A} \times (-\infty, y]\}$, $f_Y$ denotes the density of $Y$ and

$$
G_\theta(\mathcal{A}, y) = \int_{\mathcal{A}} F(y|\theta^T x) f(x)\, dx.
\tag{2.6}
$$

A low-dimensional approximation to $F_{Y|X}(y|X = x)$ is therefore $\widetilde{F}_{\hat{\theta}}(y|\hat{\theta}^T x)$, where $\widetilde{F}_\theta(y|z)$ is an estimator of $P(Y \leq y|\theta^T X = \theta^T x)$. Denoting by $\widehat{F}$ a local

linear version of $\widetilde{F}$, we define

$$(2.7) \qquad \widehat{F}_\theta(y|\theta^{\mathrm{T}}x) = \left\{\sum_{i=1}^n w_i(x,\theta)I(Y_i \leq y)\right\} \bigg/ \left\{\sum_{i=1}^n w_i(x,\theta)\right\},$$

where

$$w_i(x,\theta) = K\left\{\frac{\theta^{\mathrm{T}}(x-X_i)}{h}\right\}\left\{T^{[2]}(x,\theta) - \frac{\theta^{\mathrm{T}}(x-X_i)}{h}T^{[1]}(x,\theta)\right\},$$

$$T^{[k]}(x,\theta) = \frac{1}{nh}\sum_{i=1}^n K\left\{\frac{\theta^{\mathrm{T}}(x-X_i)}{h}\right\}\left\{\frac{\theta^{\mathrm{T}}(x-X_i)}{h}\right\}^k.$$

Our empirical, low-dimensional approximation to $F_{Y|X}(y|X=x)$ is taken to be $\widehat{F}_{\hat{\theta}}(y|\hat{\theta}^{\mathrm{T}}x)$, and is of course an estimator of $P(Y \leq y|\theta_0^{\mathrm{T}}X = \theta_0^{\mathrm{T}}x)$.

2.3. *Empirical bandwidth choice—a rule of thumb.* Two bandwidths need to be chosen: $h$ for estimating $\theta$, and $H$ for estimating $F_{\hat{\theta}}(y|\hat{\theta}^{\mathrm{T}}x)$ with $\hat{\theta}$ given. In such nonstandard problems, conventional bandwidth selection methods for nonparametric regression are either tedious to apply (as in the case of plug-in methods), or do not facilitate obvious analogies (e.g., cross-validation and its variants). Note that with $\hat{\theta}$ given, estimation of $F_{\hat{\theta}}(y|\hat{\theta}^{\mathrm{T}}x)$ has been investigated by, among others, Hall, Wolff and Yao (1999). They proposed a bootstrap method based on an approximating parametric model to determine the bandwidth, which we will adopt for estimating $H$. Furthermore, we outline a similar empirical procedure below for determining $h$.

First we fit the linear model

$$(2.8) \qquad Y_i = \beta_0 + \beta^{\mathrm{T}}X_i + \varepsilon_i.$$

Let $\check{\beta}_0$ and $\check{\beta}$ be the estimators derived by, for example, least squares, and let $\hat{\varepsilon}_1,\ldots,\hat{\varepsilon}_n$ denote the centered residuals. We shall compute a bootstrap sample $\{Y_1^*,\ldots,Y_n^*\}$ from the model

$$(2.9) \qquad Y_i^* = \check{\beta}_0 + \check{\beta}^{\mathrm{T}}X_i + \varepsilon_i^*,$$

where $\{\varepsilon_i^*\}$ denotes a conventional bootstrap resample drawn by sampling with replacement from $\{\hat{\varepsilon}_i\}$. Then the conditional distribution of $Y_i^*$, given $X_i$, depends on $X_i$ through $\check{\beta}^{\mathrm{T}}X_i$ alone. Let $\hat{\beta}^* = \hat{\beta}^*(h)$ be the estimator obtained in the same manner as $\hat{\theta}$ but with the data $(X_i,Y_i)$ replaced by their resampled counterparts $(X_i,Y_i^*)$; see Section 2.2. We choose $h$ to minimize

$$(2.10) \qquad M_1(h) = E[\|\hat{\beta}^* - \check{\beta}\|^2|\{(X_i,Y_i)\}].$$

It is important that the two bandwidths $h$ and $H$ should be different. As we shall show in Section 3, optimal performance is achieved if $h$ is of smaller order than $H$. The simulation results reported in Section 4 indicate that the bandwidths selected by the bootstrap methods discussed above produce estimators with good performance.

**3. Theory.** For simplicity we discuss only the case where the data $(X_i, Y_i)$ are independent. Analogues of our main results, Theorems 3.1 and 3.2, may be derived for dependent data, in particular for sequences of pairs $(X_i, Y_i)$ that satisfy sufficiently strong mixing conditions. The case of dependence will be explored numerically in Section 4.

Let us first define the vector of derivatives, $\dot{a}$, of a function $a$ of $\theta \in \Theta$. Let $\omega_1, \ldots, \omega_{d-1}$ be orthonormal vectors all perpendicular to $\theta$, put $\omega_{i\delta} = (1 - \delta^2)^{1/2}\theta + \delta\omega_i$ for a scalar $\delta$ and set

$$b_i = \lim_{\delta \to 0} \delta^{-1}\{a(\omega_{i\delta}) - a(\theta)\},$$

assuming the limit exists and is finite. Then

$$\dot{a}(\theta) \equiv \sum_{1 \le i \le d-1} b_i \omega_i,$$

a vector in the plane perpendicular to $\theta$. Similarly we may define the matrix, $\ddot{a}$, of second derivatives of $a$.

Let $(X, Y)$ have the distribution of a generic pair $(X_i, Y_i)$. We shall assume that

(3.1)        the density of $(X, Y)$ has four bounded derivatives, and all moments of $Y$ are finite.

The bandwidth $h$ will be permitted to vary within a range, effectively from $n^{-1/3}$ to $n^{-1/4}$; see (3.4) below. If we were confining attention to the lower end of this range, then we could reduce the smoothness assumption in (3.1) from four bounded derivatives to three derivatives plus a Hölder continuity condition. In this sense, the smoothness required by (3.1) is excessive.

Recall that if sphere radii vary in the continuum, then $\mathcal{Q}$ denotes a set of sphere centers and radii, while if there is a single, fixed radius, then $\mathcal{Q}$ is a set just of sphere centers. In either case, all spheres in $\mathcal{Q}$ are completely contained within $\mathcal{R}$; see the definition of $\mathcal{Q}$ in Section 2.2. We shall suppose that

(3.2)        $\mathcal{R}$ is an open, bounded set; the density of $X$ is bounded away from zero on $\mathcal{R}$; and the content of $\mathcal{Q}$ is nonzero.

In particular, this and (3.1) ensure that the density of the distribution of $\theta^{\mathrm{T}}X$ is bounded away from zero on the set of points $\theta^{\mathrm{T}}x$ with $x \in \mathcal{A} \subseteq \mathcal{R}$. Assumption (3.2) may therefore be viewed as the analogue of the condition, imposed in more standard problems of nonparametric regression, that the design density is bounded above zero.

Conditions (3.1) and (3.2) imply a range of smoothness properties of the marginal density $f_{\theta^{\mathrm{T}}X}$ and the conditional distribution $F(y|z) = P(Y \le y|\theta^{\mathrm{T}}X = z)$. For example, the $k_1$th derivative with respect to $\theta$, of the $k_2$th derivative with respect to $z$, of either $f_{\theta^{\mathrm{T}}X}(z)$ or $F(y|z)$, is well defined and bounded in $k_1 + k_2 \leq 4$, $y$, $\theta \in \Theta$ and $z = \theta^{\mathrm{T}}x$ for $x \in \mathcal{R}$.

Recall the definition of $G_\theta(\mathcal{A}, y)$ in (2.6), and let $\dot{G}_\theta(\mathcal{A}, y)$ and $\ddot{G}_\theta(\mathcal{A}, y)$ denote, respectively, the vector of first derivatives and the matrix of second derivatives of $G_\theta(\mathcal{A}, y)$ with respect to $\theta$, with $(\mathcal{A}, y)$ held fixed. Note that $\theta_0 = \arg\min_\theta S_0(\theta)$, where $S_0$ is defined in (2.5).

Put

$$M(\theta) = \int_{\mathcal{Q}} d\alpha \int [\dot{G}_\theta(\mathcal{A}_\alpha, y)\dot{G}_\theta(\mathcal{A}_\alpha, y)^{\mathrm{T}}$$
$$- \{F(\mathcal{A}_\alpha, y) - G_\theta(\mathcal{A}_\alpha, y)\}\ddot{G}_\theta(\mathcal{A}_\alpha, y)]f_Y(y)\,dy,$$

a $d \times d$ matrix. By assuming that

(3.3)     $\theta = \theta_0$ gives a unique global minimum of $S_0(\theta)$, and
          $\omega^{\mathrm{T}}M(\theta_0)\omega > 0$ for each nonvanishing vector $\omega \perp \theta_0$,

we require an equivalent condition that $S_0(\theta) \to S_0(\theta_0)$ at exactly the rate $\|\theta - \theta_0\|^2$ as $\theta \to \theta_0$. Of the kernel $K$ and bandwidth $h$ we shall assume that

(3.4)     $K$ is nonnegative, symmetric and compactly supported,
          and has a bounded derivative; and, for some $\varepsilon > 0$, $h = h(n)$ satisfies $h = O(n^{-\varepsilon - (1/4)})$ and $n^{-(1/3)+\varepsilon} = O(h)$ as $n \to \infty$.

The most important aspect of this assumption is that it implies $h$ should lie between $n^{-1/3}$ and $n^{-1/4}$, and so should be an order of magnitude smaller than a conventional bandwidth for estimating a univariate function by nonparametric regression. A conventional bandwidth would be of size $n^{-1/5}$.

Let $\phi_{\theta^{\mathrm{T}}X|\mathcal{A}}$ denote the density of $\theta^{\mathrm{T}}X$ conditional on $X \in \mathcal{A}$, and define

(3.5)     $\psi(\mathcal{A}, x_1, y_1, y, \theta) = \{I(y_1 \leq y) - F(y|\theta^{\mathrm{T}}x_1)\}$

$$\times \left\{I(x_1 \in \mathcal{A}) - \frac{\phi_{\theta^{\mathrm{T}}X|\mathcal{A}}(\theta^{\mathrm{T}}x_1)P(X \in \mathcal{A})}{f_{\theta^{\mathrm{T}}X}(\theta^{\mathrm{T}}x_1)}\right\}.$$

[The ratio in this definition is guaranteed well defined, since $P(X \in \mathcal{A})\phi_{\theta^{\mathrm{T}}X|\mathcal{A}} \leq f_{\theta^{\mathrm{T}}X}$.] Let $V$ denote the Gaussian $d$-vector with zero mean and covariance matrix equal to that of

$$W = \int_{\mathcal{Q}} d\alpha \left[\int \psi(\mathcal{A}_\alpha, X, Y, y, \theta)\dot{G}_{\theta_0}(\mathcal{A}_\alpha, y)f(y)\,dy\right.$$
$$\left. + \{F(\mathcal{A}_\alpha, Y) - G_{\theta_0}(\mathcal{A}_\alpha, Y)\}\dot{G}_{\theta_0}(\mathcal{A}_\alpha, Y)\right]d\alpha.$$

Let $\|\cdot\|$ denote the Euclidean metric in $d$-variate space, and recall that $\hat\theta$ is defined to be the global minimizer of $S(\theta)$ in (2.4).

THEOREM 3.1. *Assume conditions* (3.1)–(3.4). *Then* $\hat{\theta} \to \theta_0$ *with probability* 1, *and* $n^{1/2} M(\theta_0)(\hat{\theta} - \theta_0)$ *converges in distribution to* $V$ *as* $n \to \infty$.

To appreciate the implications of this result, let $\hat{\theta}^{\perp}$ denote the projection of $\hat{\theta}$ into the plane $\Pi^{\perp}$ that is perpendicular to $\theta_0$. (Equivalently, $\hat{\theta}^{\perp}$ is the projection of $\hat{\theta} - \theta_0$ into $\Pi^{\perp}$.) The first part of Theorem 3.1 implies that $\|\hat{\theta} - \theta_0\| \to 0$ with probability 1, from which it follows (since $\hat{\theta}$ and $\theta_0$ are both unit vectors) that

$$(3.6) \qquad \hat{\theta} - \theta_0 = \hat{\theta}^{\perp} + o(\|\hat{\theta} - \theta_0\|)$$

with probability 1. That is, in first-order asymptotic terms, $\hat{\theta} - \theta_0$ is completely describable through the projection of this vector into the plane perpendicular to $\theta_0$.

Note that, by definition of differentiation with respect to $\theta$, the vector $\dot{G}_{\theta}$ is perpendicular to $\theta$. It therefore follows from the definition of $V$ that, with probability 1, $V$ lies completely in $\Pi^{\perp}$. Observe too that, in view of (3.3), there is a generalized inverse of $M_0 = M(\theta_0)$ (call it $M_0^{-}$) that is well defined in $\Pi^{\perp}$. It has the property that

$$M_0 M_0^{-} v = M_0^{-} M_0 v = v \qquad \text{for all } v \in \Pi^{\perp}.$$

These results, Theorem 3.1 and (3.6) imply that $n^{1/2}(\hat{\theta} - \theta)$ converges in distribution to $M_0^{-} V$.

Of course, our main purpose in computing $\hat{\theta}$ is so it can be used in a conditional distribution estimator, such as $\hat{F}_{\theta}$ introduced in (2.7). Theorem 3.2 below shows that the root-$n$ consistency achieved by the estimator $\hat{\theta}$ makes that quantity so accurate that, from the viewpoint of first-order performance, the estimator $\hat{F}_{\hat{\theta}}(y|\hat{\theta}^{\mathrm{T}} x)$ is equivalent to its counterpart which would be employed if the value of $\theta_0$ were known. This result has analogues for general choice of the bandwidth used for $\hat{F}_{\theta}$; they describe a range of circumstances where the leading bias and variance terms do not include the effect of estimating $\theta$. However, for the sake of simplicity and brevity we shall treat only the optimal size of bandwidth.

The latter size is $n^{-1/5}$, and when that is employed, $\hat{F}_{\theta_0}(y|\theta_0^{\mathrm{T}} x)$ converges to its limit at rate $n^{-2/5}$. We shall show in Theorem 3.2 that the difference between $\hat{F}_{\hat{\theta}}(y|\hat{\theta}^{\mathrm{T}} x)$ and $\hat{F}_{\theta_0}(y|\theta_0^{\mathrm{T}} x)$ is then of strictly smaller order than $n^{-2/5}$.

These considerations motivate the following assumption:

$$(3.7) \qquad \begin{array}{l} \text{the bandwidth } H \text{ used to construct } \hat{F}_{\theta} \text{ has the property} \\ \text{that } n^{1/5} H \text{ is bounded away from zero and infinity as } n \to \\ \infty; \text{ and the kernel is nonnegative, symmetric, compactly} \\ \text{supported and has a bounded derivative.} \end{array}$$

*Note that $H$ and $h$ are of different orders, the former being of size $n^{-1/5}$ and the latter of smaller order.* We shall reduce the stringency of (3.1), assuming instead that

(3.8)   the density of $(X, Y)$ has two continuous derivatives, and all moments of $Y$ are finite.

As the following theorem shows, we do not need the full force of the result that $\hat{\theta} - \theta_0 = O_p(n^{-1/2})$; the convergence rate $o_p(n^{-2/5})$ suffices.

THEOREM 3.2.  *Assume* (3.2), (3.7), (3.8), *that $x \in \mathcal{R}$, and that $\hat{\theta} - \theta_0 = o_p(n^{-2/5})$ as $n \to \infty$. Then for each $y$,*

$$\widehat{F}_{\hat{\theta}}(y|\hat{\theta}^{\mathrm{T}}x) = \widehat{F}_{\theta_0}(y|\theta_0^{\mathrm{T}}x) + o_p(n^{-2/5}).$$

It follows from the asymptotic normality of local linear regression estimation [see, e.g., Theorem 1 of Fan, Heckman and Wand (1995), and Remark 4 of Hall, Wolff and Yao (1999)] that the estimator $\widehat{F}_{\theta_0}(y|\theta_0^{\mathrm{T}}x)$ is asymptotically normally distributed with convergence rate $n^{-2/5}$. By Theorem 3.2 above, $\widehat{F}_{\hat{\theta}}(y|\hat{\theta}^{\mathrm{T}}x)$ and $\widehat{F}_{\theta_0}(y|\theta_0^{\mathrm{T}}x)$ have the same asymptotic distribution.

**4. Numerical properties.**  We approximate the integral in (2.4) by a series,

(4.1)   $$S(\theta) = \frac{1}{B}\sum_{i=1}^{B} S(\theta, \mathcal{A}_i),$$

where the $\mathcal{A}_i$'s are spheres of radius $r$ contained within $\mathcal{R}$. In practice one would select a value of $B$ that permitted the calculations to be completed within a reasonable time, and compute estimates for that value as well as for substantially smaller ones, say half and three-quarters of the initial $B$. Provided there was little variation in the results, the larger $B$ would be appropriate. The results reported in this section show that choice of $B$ has little effect on final results.

In the numerical examples below we searched for $\theta$ (with $h$ fixed) using the downhill simplex method; see Section 10.4 of Press, Teukolsky, Vetterling and Flannery (1992). Using the Epanechnikov kernel, the bandwidths were sought among values $h_i = 0.1 \times 1.2^{i-1}$ for $i = 1, \ldots, 15$, based on the bootstrap methods outlined in Section 2.3. We used sample sizes $n = 200$ and $400$. Each setting was replicated 50 times. Throughout Examples 1 and 2 below we took $X_{ij}$ and $\varepsilon_i$ to be totally independent $N(0, 1)$ random variables.

EXAMPLE 1.  Here we consider the model

$$Y_i = \theta_1 X_{i1} + \theta_2 X_{i2} + \theta_3 X_{i3} + \theta_4 X_{i4} + \varepsilon_i,$$

where $\theta^{\mathrm{T}} \equiv (\theta_1, \ldots, \theta_4) = (1, 2, 0, 3)/\sqrt{14}$. Thus, the conditional distribution of $Y$, given $X \equiv (X_1, \ldots, X_4)^{\mathrm{T}}$, is $\mathrm{N}(\theta^{\mathrm{T}} X, 1)$. We let the radius be $r = 1$, and sphere centers be points $(x_1, x_2, x_3, x_4)$, where each $x_j$ ranged over either five or seven grid points between $-1.5$ and $1.5$, with spacing $0.75$ or $0.5$, respectively, resulting in $B = 625$ or $B = 2401$.

Figure 1(a)–(c) presents boxplots of the inner product $\theta^{\mathrm{T}}\hat{\theta}$, where, respectively, the bandwidth $h$ was computed by minimizing (2.10), or taken equal to the latter value multiplied by 1.5 or 0.7. Since both $\theta$ and $\hat{\theta}$ are unit vectors, $\theta^{\mathrm{T}}\hat{\theta} = 1$ if and only if $\theta = \hat{\theta}$. We see from Figure 1(a)–(c) that the estimates of $\theta$ become steadily more accurate as sample size increases. Moreover, the algorithm is largely insensitive to the bandwidths used in the search; the estimates of $\theta$ with the three different bandwidths differ only a little. Furthermore, the algorithm is also insensitive to the value of $B$.

Figure 1(d) and (e) displays boxplots of the bandwidths $h$, obtained by minimizing $M_1$ in (2.10), and $H$, defined by the method of Hall, Wolff and Yao (1999). As expected, empirical bandwidth is a decreasing function of sample size. Note too that selected $h$'s are in general noticeably smaller than the chosen $H$, which is in agreement with the asymptotic orders of $h$ and $H$.

We also calculated values of the local linear estimator defined in (2.7) with bandwidth $H$. Figure 1(f) gives average absolute errors, computed using a regular grid (with adjacent points distant 0.05 apart) in the $(\theta^{\mathrm{T}} X, Y)$-plane. For the sake of comparison we also report the errors for the estimators based on the true $\theta$. Clearly, accuracy increases with sample size, and estimators based on $\hat{\theta}$ are less accurate than those based on the true $\theta$. However, the deficit due to errors in estimating $\theta$ is not great when $n = 200$, and is negligible when $n = 400$. Choice of radius $r$ is not critical either; results with $r = 0.5$ and 1.5 are similar to those for $r = 1$, and therefore are not reported here.

EXAMPLE 2. Next we consider the model

$$Y_i = \tfrac{1}{2}(\sin X_{i1} + \sin X_{i2} + \sin X_{i3} + \sin X_{i4}) + \varepsilon_i.$$

Now the conditional distribution of $Y$ given $X = (X_1, \ldots, X_4)^{\mathrm{T}}$ no longer depends on a linear combination of $X$. The true value of $\theta$ is $(0.5, 0.5, 0.5, 0.5)^{\mathrm{T}}$; note the symmetry of the model. We selected the spheres in the same way as in Example 1. The numerical results are presented in Figure 2, which displays a similar pattern to Figure 1 although the estimates in general are not as accurate as in Example 1. This is due to the fact that we were estimating the least-squares approximation, in the sense of minimizing (4.1), of the conditional distribution of $Y$ given $X$, rather than the conditional distribution itself. Figure 2(a)–(c) shows that the estimation for $\theta$ is still accurate, even for the sample size $n = 200$, and is steadily improved when $n$ is increased to 400.

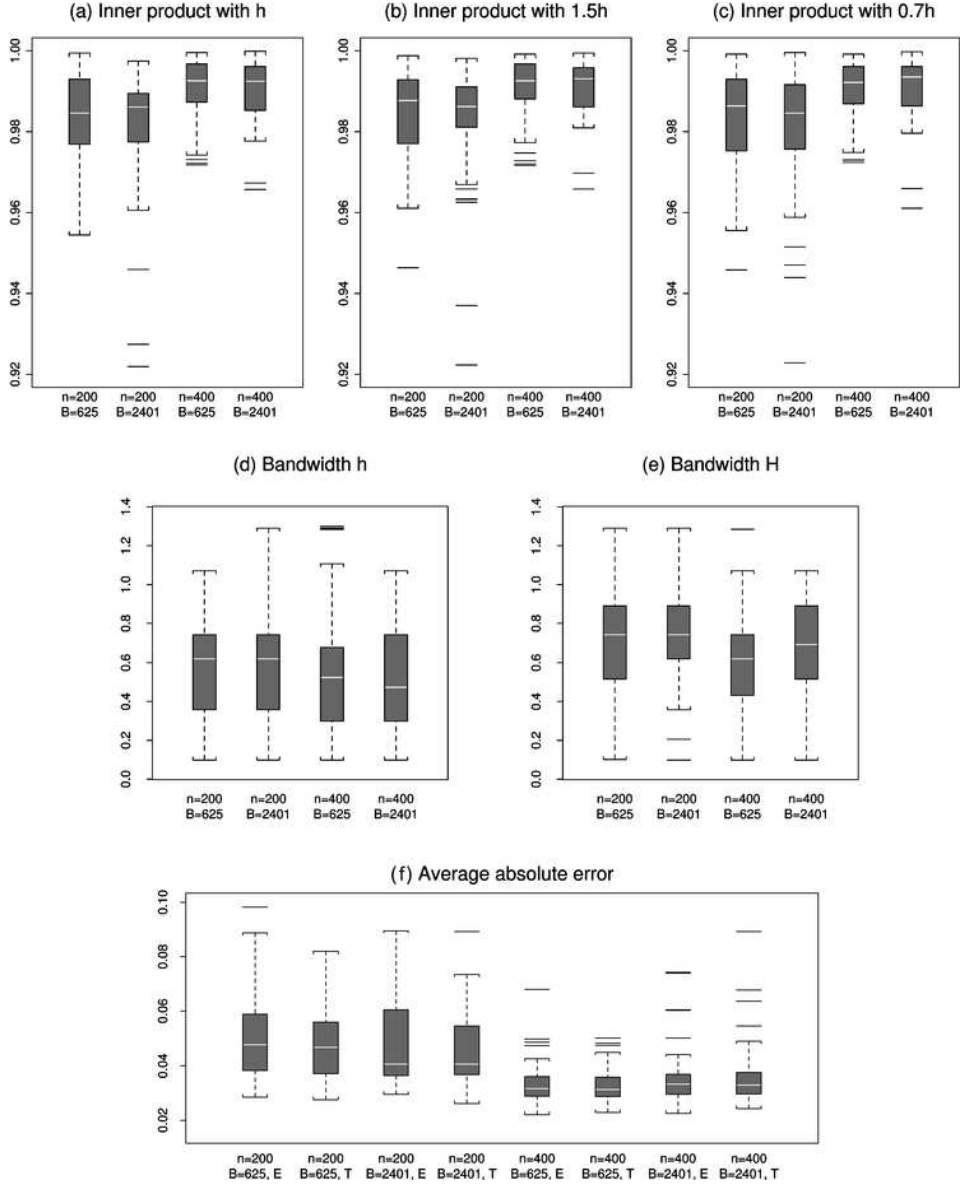

Fig. 1.  *Simulation results for Example* 1. *Boxplots of the inner product* $\hat{\theta}^{\mathrm{T}}\theta$, *with bandwidth* $h$ *taken equal to* (a) $\hat{h}$, (b) $1.5\hat{h}$, *and* (c) $0.7\hat{h}$; *and of* (d) $\hat{h}$, (e) $\hat{H}$, *and* (f) *average absolute errors of estimated conditional distribution of* $Y$ *given* $\theta^{\mathrm{T}}X$ *with either* $\theta = \hat{\theta}$ *(denoted by "E") or* $\theta$ *equal to its true value (denoted by "T").*

EXAMPLE 3.  Finally we illustrate our method with $\{Y_t,\ 1 \leq t \leq 176\}$ the quarterly growth rates of US real GNP between February 1947 and

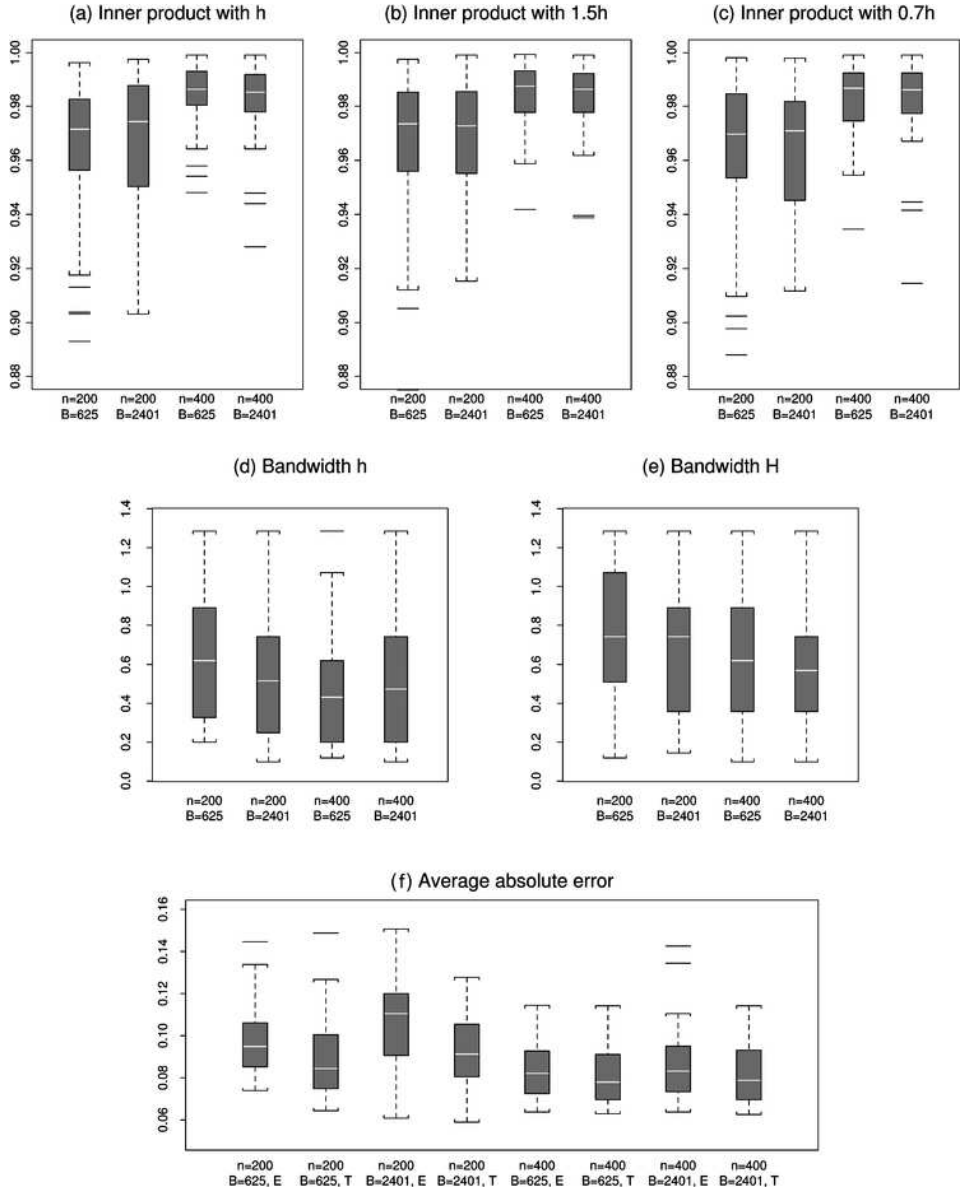

FIG. 2. *Simulation results for Example 2. Panels show the same information as in Figure 1.*

January 1991. The data series is plotted in Figure 3. This dataset has been analyzed by, for example, Tiao and Tsay (1994). Let $X_t = (Y_{t-1}, Y_{t-2})^{\mathrm{T}}$. We estimated the value of $\theta = (\theta_1, \theta_2)^{\mathrm{T}}$ for which the conditional distribution of $Y_t$, given $\theta^{\mathrm{T}} X_t$, was the best approximation for the conditional distribution of $Y_t$ given $X_t$, in the sense that $S(\theta)$, defined at (4.1), was minimized. We first standardized the data $X_t$. Sphere centers were taken to be the points $X_t$ (so that $B = n$), with radius $r = 1$. The resulting estimate is $\hat{\theta} = (0.580, -0.815)^\mathrm{T}$.

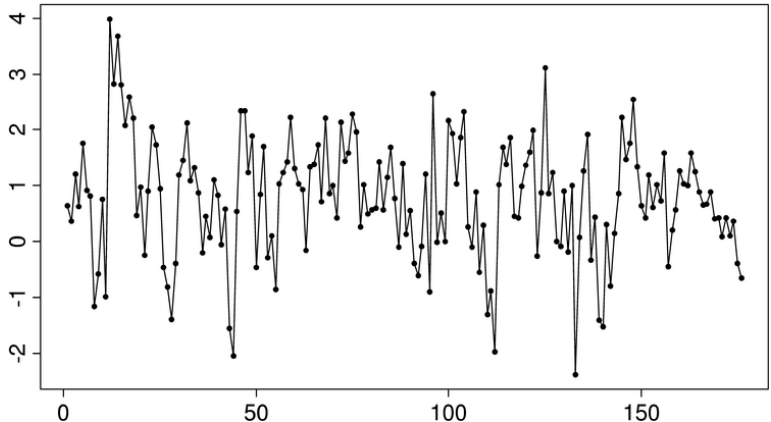

Fig. 3.   *Prediction intervals for US quarterly GNP growth $X_t$ based on, respectively, three different predictors $0.580X_{t-1} - 0.815X_{t-2}$, $X_{t-1}$ and $(X_{t-1}, X_{t-2})$.*

Once $\hat{\theta}$ was obtained we constructed the adjusted Nadaraya–Watson estimator $\widehat{F}(\cdot|z)$ [see Hall, Wolff and Yao (1999)] of the conditional distribution of $Y_t$, given $\hat{\theta}^\mathrm{T} X_t = z$. The resulting quantile prediction interval is $[\widehat{F}^{-1}(\frac{1}{2}\alpha|z), \widehat{F}^{-1}(1 - \frac{1}{2}\alpha|z)]$, for $\alpha \in (0,1)$. To check on performance we used the first 166 data points to estimate $\hat{\theta}$ and $\widehat{F}(\cdot|z)$, and employed the last ten data points to validate the predicted values. Results with $\alpha = 0.1$ are reported in Table 1. Note that with $\hat{\theta} = (0.580, -0.815)^\mathrm{T}$, the predictor is $0.580Y_{t-1} - 0.815Y_{t-2}$. For comparison we also report prediction intervals using a single predictor $Y_{t-1}$, and a two-dimensional predictor $(Y_{t-1}, Y_{t-2})$.

All the intervals in the table contain the corresponding true values. Prediction intervals based on two predictors $Y_{t-1}$ and $Y_{t-2}$ are more accurate, in general, than those based on a single predictor $Y_{t-1}$, since the average length of the prediction intervals is reduced from 3.51 to 3.21. It is interesting to see that the average length of the prediction intervals based on the selected single predictor $0.580Y_{t-1} - 0.815Y_{t-2}$ is 3.22, which is almost the same as that based on $(Y_{t-1}, Y_{t-2})$. Note too that our method does not use multivariate smoothing techniques, which are susceptible to the "curse of dimensionality." Predictions based on $d = 3$ and 4 did not lead to significant improvements, and therefore are omitted. The absence of improvement is in agreement with results of Tiao and Tsay (1994), who proposed nonlinear, second-order autoregressive models for this dataset.

## 5. Outlines of technical arguments.

OUTLINE PROOF OF THEOREM 3.1.   Our argument has two main stages, showing, respectively, that

$$(5.1) \qquad \|\hat{\theta} - \theta_0\| = O_p(n^{\varepsilon - (1/2)}) \qquad \text{for each } \varepsilon > 0,$$

$$(5.2) \qquad \begin{aligned} S(\theta) = {} & T + (\theta - \theta_0)^{\mathrm{T}} M_0 (\theta - \theta_0) - 2(\theta - \theta_0)^{\mathrm{T}} (V_1 + V_2) \\ & + o_p(\|\theta - \theta_0\|^2) + O_p\Big( \sum \|\theta - \theta_0\|^u n^{-v-\zeta} \Big), \end{aligned}$$

where $T$ does not depend on $\theta$, $\zeta > 0$ is fixed,

$$V_1 = \int d\widehat{F}(y) \int_{\mathcal{R}_r} \dot{G}_{\theta_0}(\mathcal{A}_\alpha, y) \xi_n(\mathcal{A}_\alpha, y, \theta_0) \, d\alpha,$$

$$V_2 = \int d\widehat{F}(y) \int_{\mathcal{R}_r} D_{\theta_0}(\alpha, y) \dot{G}_{\theta_0}(\mathcal{A}_\alpha, y) \, d\alpha,$$

$$\xi_n(\mathcal{A}, y, \theta) = \frac{1}{n} \sum_{i=1}^n [\psi(\mathcal{A}, X_i, Y_i, y, \theta) - E\{\psi(\mathcal{A}, X, Y, y, \theta)\}],$$

$\psi$ is as in (3.5), and $O_p(\sum \|\theta - \theta_0\|^u n^{-v-\zeta})$ denotes a quantity which uniformly in $\theta$ is of order no more than that of the sum of $\|\theta - \theta_0\|^u n^{-v-\zeta}$ over a fixed, finite set of pairs $(u, v)$, where in each case, $u, v \geq 0$ and $\frac{1}{2} u + v \geq 1$.

To give an appreciation of the origin of the terms which make up the $O_p(\cdots)$ remainder in (5.2), we note that the contributions to the remainder come from different steps in a Taylor expansion of $S(\theta)$. In particular, terms of the following orders arise in that way:

$$(5.3) \qquad \begin{aligned} & \|\theta - \theta_0\| h^2, \qquad \|\theta - \theta_0\| (nh^{3/2})^{-1} n^\varepsilon, \qquad \|\theta - \theta_0\| n^{-t-(1/2)}, \\ & \|\theta - \theta_0\|^2 (nh^3)^{-1}, \qquad (nh^3)^{\varepsilon - 2}, \qquad n^{-t-1}, \end{aligned}$$

TABLE 1

| True value | $0.580X_{t-1} - 0.815X_{t-2}$ | $X_{t-1}$ | $(X_{t-1}, X_{t-2})$ |
|---|---|---|---|
| 0.67 | $[-0.99, 2.32]$ | $[-0.99, 2.32]$ | $[-0.62, 3.11]$ |
| 0.89 | $[-0.91, 2.32]$ | $[-0.88, 2.34]$ | $[-0.59, 2.28]$ |
| 0.40 | $[-0.99, 2.20]$ | $[-1.56, 2.54]$ | $[-0.86, 2.34]$ |
| 0.43 | $[-0.91, 2.34]$ | $[-0.99, 2.32]$ | $[-0.62, 3.11]$ |
| 0.09 | $[-0.91, 2.28]$ | $[-0.88, 2.34]$ | $[-0.59, 2.21]$ |
| 0.42 | $[-0.99, 2.20]$ | $[-1.56, 2.54]$ | $[-1.17, 2.34]$ |
| 0.11 | $[-0.88, 2.32]$ | $[-0.99, 2.32]$ | $[-0.62, 2.32]$ |
| 0.36 | $[-0.91, 2.34]$ | $[-0.88, 2.34]$ | $[-0.59, 2.12]$ |
| $-0.40$ | $[-0.99, 2.34]$ | $[-1.56, 2.54]$ | $[-0.86, 2.54]$ |
| $-0.65$ | $[-0.81, 2.32]$ | $[-0.91, 2.32]$ | $[-0.91, 2.32]$ |
| Average length | 3.22 | 3.51 | 3.21 |

where in each case the bound is valid for all $\varepsilon > 0$ and some $t > 0$. Noting that, by (3.4), $n^{\zeta_1 - (1/3)} \leq h \leq n^{-\zeta_2 - (1/4)}$ for constants $\zeta_1, \zeta_2 > 0$, and using the upper of these bounds when $h$ appears with a positive exponent in (5.3), and the lower when $h$ appears with a negative exponent, we see that each of the quantities in (5.3) may be written as $\|\theta - \theta_0\|^u n^{-v-\zeta}$ for some $\zeta > 0$ and some $(u, v)$ such that $\frac{1}{2}u + v \geq 1$.

More detailed proofs of (5.1) and (5.2) can be found in Hall and Yao (2002). To illustrate the use of the regularity conditions (3.1)–(3.4), we mention that (3.1) is employed to guarantee adequate smoothness of $F$ when Taylor-expanding $F(Y_j | \theta^T x)$ and related functions; that (3.1) and (3.2) together ensure that the effective design density is bounded away from zero, which allows us to deal with the denominator of $\widehat{F}_{-i,-j}(Y_j | \theta^T X_i)$ via a stochastic Taylor expansion; that (3.3) guarantees that the minimum of $S(\theta)$ is attained in the usual quadratic way, or equivalently that the matrix $M(\theta_0)$ is of full rank in the $(d-1)$-dimensional space of vectors perpendicular to $\theta_0$; and that one of the applications of (3.4) was described in the previous paragraph.

Taking $\theta = \hat{\theta}$ in (5.2), and noting (5.37), we see that the remainder term in (5.2) may be written as $O_p(\sum n^{-(u/2)-v-\zeta})$. Since $\frac{1}{2}u + v \geq 1$ for each pair $(u, v)$ contributing to the series, and since $\zeta > 0$, then this $O_p(\cdots)$ remainder equals $o_p(n^{-1})$. Theorem 3.1 follows from this form of (5.2), and from the fact that $n^{1/2}(V_1 + V_2)$ converges in distribution to $V$, the latter defined a little before the statement of the theorem.   □

OUTLINE PROOF OF THEOREM 3.2.   Let $\Theta_n$ denote the set of all $\theta \in \Theta$ that satisfy $\|\theta - \theta_0\| \leq \delta(n) n^{-2/5}$, where $\delta(n) \downarrow 0$ as $n \to \infty$. The theorem follows from the following result.

LEMMA.   Assume (3.2), (3.7), (3.8) and that $x \in \mathcal{R}$. Then for each $y$

$$\sup_{\theta \in \Theta_n} |\widehat{F}_\theta(y | \theta^T x) - \widehat{F}_{\theta_0}(y | \theta_0^T x)| = o_p(n^{-2/5}).$$

We outline the proof of the lemma. Treat $\widehat{F}_\theta$ as the ratio expressed in (2.7), although multiply top and bottom there by $(nh)^{-1}$ [here $(nH)^{-1}$, since we take the bandwidth to be $H$] in order to ensure that neither the numerator nor the denominator converges to zero or diverges to infinity. The numerator and denominator are now each in the form $T_1 T_2 - T_3 T_4$, where each $T_j$ is linear in functions of the data $X_i$ and has a proper limit as $n$ diverges. Additively decompose each $T_j$ into its expected value (or mean), and the difference between it and its mean. Each mean is of course purely deterministic. In the remainder of this section we shall outline the technique, starting from this decomposition, for treating $T_1$ and $T_2$; a similar argument may be given in the case of $T_3$ or $T_4$.

The expected value of $T_1$ or $T_2$ may be written as its "$H \to 0$ limit," plus a term that equals $H^2$ multiplied by a function of $\theta$, plus a remainder that equals $o(H^2)$ uniformly in $\theta$. The "$H \to 0$ limit," evaluated at $\theta$, equals the same quantity evaluated at $\theta_0$ rather than at $\theta$, plus a remainder of order $O\{\delta(n)n^{-2/5}\} = o(n^{-2/5})$, uniformly in $\theta \in \Theta_n$; and similarly, the coefficients of $H^2$ (for $\theta$ and $\theta_0$, resp.) are identical, up to a term that converges to 0 uniformly in $\theta \in \Theta_n$ as $n \to \infty$. These arguments require only Taylor expansion, and prove that the mean of each of the $T_j$'s equals its counterpart when $\theta$ is replaced by $\theta_0$, plus terms that are of size $o(n^{-2/5})$ uniformly in $\theta \in \Theta_n$. A longer argument [see Hall and Yao (2002)] can be used to show that the same property is enjoyed by each $T_j - E(T_j)$, not just by each $E(T_j)$. The theorem follows from these properties. $\square$

**Acknowledgments.** We are grateful to an Editor and two reviewers for helpful comments.

CENTRE FOR MATHEMATICS
  AND ITS APPLICATIONS
AUSTRALIAN NATIONAL UNIVERSITY
CANBERRA, ACT 0200
AUSTRALIA
E-MAIL: Peter.Hall@maths.anu.edu.au

DEPARTMENT OF STATISTICS
LONDON SCHOOL OF ECONOMICS
HOUGHTON STREET
LONDON WC2A 2AE
UNITED KINGDOM
E-MAIL: Q.Yao@lse.ac.uk